\documentclass[12pt]{article}
\usepackage{amsthm,amsfonts,amssymb,amscd}

\begin{document}

\begin{center}
\large \bf Birationally rigid finite covers \\
of the projective space
\end{center}\vspace{0.5cm}

\centerline{A.V.Pukhlikov}\vspace{0.5cm}

\parshape=1
3cm 10cm \noindent {\small \quad\quad\quad \quad\quad\quad\quad
\quad\quad\quad {\bf }\newline In this paper we prove birational superrigidity of finite covers of degree $d$ of the $M$-dimensional projective space of index 1, where $d\geqslant 5$ and $M\geqslant 10$, with at most quadratic singularities of rank $\geqslant 7$, satisfying certain regularity conditions. Up to now, only cyclic covers were studied in this respect. The set of varieties with worse singularities or not satisfying the regularity conditions is of codimension $\geqslant\frac12(M-4)(M-5)+1$ in the natural parameter space of the family.

Bibliography: 14 titles.} \vspace{1cm}

\vspace{1cm}

\noindent 14E05, 14E07\vspace{1cm}

\section*{Introduction}

{\bf 0.1. Statement of the main result.}  Let us fix integers $d\geqslant 5$ and $l\geqslant 2$, where $(d,l)\neq(5,2)$. Set $M=(d-1)l$, so that $M\geqslant 10$. In the present paper we study
$d$-sheeted covers of the complex projective space
${\mathbb P}={\mathbb P}^M$ with at most quadratic singularities of rank $\geqslant 7$, which are Fano varieties of index 1. Such covers have a convenient presentation: let
$$
\overline{{\mathbb P}}={\mathbb P}(\underbrace{1,\dots,1}_{M+1},l)={\mathbb P}(1^{M+1},l)
$$
be the weighted projective space with homogeneous coordinates $x_0,\dots,x_M$, $\xi$, where $x_i$ are of weight 1 and $\xi$ is of weight $l$. Furthermore, let
$$
F(x_*,\xi)=\xi^d+A_1(x_*)\xi^{d-1}+\dots+A_d(x_*)
$$
be a (quasi)homogeneous polynomial of degree $dl$, that is, $A_i(x_0,\dots,x_M)$ is a homogeneous polynomial of degree $il$ for $i=1,\dots,d$. The space
$$
{\cal F}=\prod^d_{i=1}H^0({\mathbb P},{\cal O}_{\mathbb P}(il))
$$
parametrizes all such polynomials. If the hypersurface
$$
V=\{F=0\}\subset\overline{\mathbb P}
$$
has at most quadratic singularities of rank $\geqslant 7$, then the set $\mathop{\rm Sing}V$ of singular points is of codimension $\geqslant 6$ in $V$, so that by Grothendieck's theorem \cite{CL} the variety $V$ is factorial. Since the property to have at most quadratic singularities of rank $\geqslant r$ is stable with respect to blow ups (see \cite[Subsection 3.1]{Pukh15a}), the singularities of the variety  $V$ are terminal. Now
$$
\mathop{\rm Pic}V={\mathbb Z}H,\quad K_V=-H,
$$
where $H$ is the class of a hyperplane section, so that $V$ is a primitive Fano variety of index 1.\vspace{0.1cm}

Let $o^*=(0:\dots:0:1)=(0^{M+1}:1)\in\overline{\mathbb P}$ be the unique singular point of the weighted projective space $\overline{\mathbb P}$. Obviously, $o^*\not\in V$. Consider the projection
$$
\pi_{\mathbb P}\,\colon\overline{\mathbb P}\backslash\{o^*\}\to{\mathbb P},
$$
$\pi_{\mathbb P}((x_0:\dots:x_M:\xi))=(x_0:\dots:x_M)$ ``from the point $o^*$''. Obviously,
$$
\pi=\pi_{\mathbb P}|_V\colon V\to{\mathbb P}
$$
is a $d$-sheeted ramified cover of the projective space. (In particular, $H$ is the $\pi$-pull back of the class of a hyperplane in ${\mathbb P}$ onto
$V$.)\vspace{0.1cm}

Now let us state the main result of the paper. We identify a polynomial
$F\in{\cal F}$ and the corresponding closed set $\{F=0\}$, which enables us to write $V\in{\cal F}$.\vspace{0.1cm}

{\bf Theorem 0.1.} {\it There is a Zariski open subset ${\cal F}_{\rm reg}\subset{\cal F}$, such that:\vspace{0.1cm}

{\rm (i)} every $V\in{\cal F}_{\rm reg}$ is a factorial Fano variety of index 1 with terminal singularities,\vspace{0.1cm}

{\rm (ii)} the inequality
$$
\mathop{\rm codim}(({\cal F}\backslash{\cal F}_{\rm reg})\subset{\cal F})\geqslant\frac12(M-4)(M-5)+1,
$$
holds,\vspace{0.1cm}

{\rm (iii)} every variety $V\in{\cal F}_{\rm reg}$ is birationally superrigid.}\vspace{0.1cm}

{\bf Corollary 0.1.} {\it Let $V\in{\cal F}_{\rm reg}$. The following claims are true.\vspace{0.1cm}

{\rm (i)} Every birational map $\chi\colon V\dashrightarrow V'$ onto a Fano variety $V'$ with ${\mathbb Q}$-factorial terminal singularities and Picard number 1 is a (biregular) isomorphism.\vspace{0.1cm}

{\rm (ii)} There is no rational dominant map $V\dashrightarrow S$ onto a positive-dimensional variety $S$, the fibre of which has negative Kodaira dimension. Therefore, on $V$ there are no structures of a rationally connected fibre space and Mori fibre space over a positive-dimensional base. In particular, $V$ has no structures of a conic bundle and $V$ is non-rational.\vspace{0.1cm}

{\rm (iii)} The groups of birational and biregular automorphisms of the variety $V$ coincide:}
$$
\mathop{\rm Bir}V=\mathop{\rm Aut}V.
$$

{\bf Proof.} All these claims are the standard implications of the birational superrigidity, see \cite[Chapter 2, Section 1]{Pukh13a}.\vspace{0.1cm}

Note that every automorphism of the variety $V\in{\cal F}_{\rm reg}$ is induced by an automorphism of the ambient weighted projective space $\overline{{\mathbb P}}$.\vspace{0.3cm}


{\bf 0.2. The structure of the paper.}  In Subsection 0.3 we give a list of known results about birational superrigidity of finite covers of index 1. All previous results were about {\it cyclic} covers (for cyclic covers the standard procedure of constructing hypertangent divisors worked well \cite[Chapter 3]{Pukh13a}, whereas in the case of an arbitrary cover this is not the case).\vspace{0.1cm}

In \S 1 we give a precise definition of the set ${\cal F}_{\rm reg}$. This definition includes several conditions, one of which is the condition to have at most quadratic singular points of rank $\geqslant 7$. Besides, we need at every point $o\in V$ a certain {\it regularity condition}, which is similar, but not identical to the usual regularity conditions on which the technique of hypertangent divisors is based.\vspace{0.1cm}

In \S 2 we prove part (ii) of Theorem 0.1. As usual (see \cite[Chapter 3]{Pukh13a}), we estimate the codimension of the set of hypersurfaces $V$, containing a fixed point $o$ and not regular at that point. After that it is not difficult to globalize the estimate for the codimension.\vspace{0.1cm}

In \S 3 we prove the birational superrigidity of regular hypersurfaces $V$. Assuming that the claim (iii) of Theorem 0.1 is not true, we obtain the existence of a mobile linear system $\Sigma\subset|nH|$ with a maximal singularity. In order to prove the birational superrigidity, we have to exclude all possible types of maximal singularities. The main technical ingredients are the $8n^2$-inequality for a non-singular point $o$ of the hypersurface $V$ and the recently discovered generalized $4n^2$-inequality for a complete intersection singularity, and, of course, the hypertangent divisors. It is not possible to apply the well known technique of hypertangent divisors directly to non-cyclic covers; the essence of this paper is precisely to modify that technique, applying it not to the variety  $V$ it self, but to its intersection with another hypersurface in the weighted projective space.\vspace{0.3cm}


{\bf 0.3. Historical remarks and acknowledgements.} The double covers of the projective space of index 1 are the ideal objects for the theory of birational rigidity, due to their low degree. Soon after the classical paper \cite{IM}, the birational superrigidity of non-singular double covers of the projective space ${\mathbb P}^3$, branched over a sextic (``sextic double solids'') was shown in \cite{I80}. In the arbitrary dimension the birational superrigidity of non-singular double spaces of index 1 was proved in \cite{Pukh89a}, and with certain types of singularities in the papers \cite{Pukh97b,Ch04b,Ch07a,ChPark10,Mullany}. The cyclic covers of arbitrary degree were studied in \cite{Pukh09a} in a much more general context, triple cyclic covers with double points were considered in \cite{Ch04c}. However, up to the present paper, non-cyclic covers were never studied. The reason was that the technique of hypertangent divisors can not be applied directly to non-cyclic covers, because in the weighted projective presentation there is a coordinate of weight $\geqslant 2$. The aim pf the present paper is to overcome this issue.\vspace{0.1cm}

The author thanks The Leverhulme Trust for the support of this work (Research Project Grant RPG-2016-279).\vspace{0.1cm}

The author is also grateful to the colleagues in the Divisions of Algebraic Geometry and Algebra of Steklov Mathematical Institute for the interest to his work, and also to the colleagues-algebraic geometers at the University of Liverpool for general support.

 \section{Regular varieties}

In this section we carry out some preliminary work, which we need to study finite covers of the projective space. In Subsection 1.1 we consider the systems of affine and homogeneous coordinates on $\overline{{\mathbb P}}$ and hypersurfaces in $\overline{{\mathbb P}}$ that ``replace'' hyperplanes of the ordinary projective space. In Subsection 1.2 we consider in more details the local equation of the hypersurface $V\in{\cal F}$ with respect to the affine coordinates on a suitable open subset in $\overline{{\mathbb P}}$. On that basis, in Subsection 1.3 we state the regularity conditions, defining the subset $\cal F_{\rm reg}\subset{\cal F}$. The claim (i) of Theorem 0.1 follows immediately from the statement of these conditions.\vspace{0.3cm}

{\bf 1.1. Preliminary remarks.} Let
$$
f(z_*)=q_{\mu}(z_*)+\dots+q_N(z_*)
$$
be a polynomial in the variables $z_1,\dots,z_M$, decomposed into homogeneous components $q_i$ of degree $i\geqslant 1$ (so that $f(0,\dots,0)=0$). Set
$o=(0,\dots,0)\in{\mathbb A}^M_z$.\vspace{0.1cm}

{\bf Definition 1.1.} The affine hypersurface $\{f=0\}$ is $k$-{\it regular at the point} $o$, where $\mu\leqslant k\leqslant N$, if the homogeneous polynomials
$$
q_{\mu},\dots,q_k
$$
form a regular sequence in ${\cal O}_{o,{\mathbb A}^M}$.\vspace{0.1cm}

Obviously, the condition of $k$-regularity at the point $o$ means that the system of equations
$$
q_{\mu}=\dots =q_k=0
$$
defines a closed subset (a cone) of codimension $k-\mu+1$ in ${\mathbb A}^M$. This conditions is meaningful only for $k-\mu+1\leqslant M$.\vspace{0.1cm}

Now let us consider hypersurfaces in the weighted projective space
$\overline{\mathbb P}$.\vspace{0.1cm}

{\bf Proposition 1.1.} {\it For every homogeneous polynomial $\gamma(x_0,\dots,x_M)$ of degree $l$ the equation $\xi=\gamma(x_*)$ defines a hypersurface $R_{\gamma}\subset\overline{\mathbb P}$ that does not contain the point $o^*=(0^{M+1}:1)$. The projection $\pi_{\mathbb P}|_{R_{\gamma}}$ is an isomorphism of $R_{\gamma}$ and} ${\mathbb P}={\mathbb P}^M$.\vspace{0.1cm}

{\bf Proof:} This is obvious.\vspace{0.1cm}

In the affine chart $\{x_0\neq 0\}\subset\overline{\mathbb P}$ with the natural affine coordinates $z_i=x_i\slash x_0$ and $y=\xi\slash x^l_0$ the projection $\pi_{\mathbb P}$ takes the form of the usual projection
$$
{\mathbb A}^{M+1}_{z_1,\dots,z_M,y}\to{\mathbb A}^M_{z_1,\dots,z_M},
$$
$$
(z_1,\dots,z_M,y)\mapsto(z_1,\dots,z_M),
$$
where ${\mathbb A}^M_z$ is the affine chart $\{x_0\neq 0\}$ in ${\mathbb P}$. Obviously, the affine hypersurface $R_{\gamma}\cap\{x_0\neq 0\}$ is given by the equation
$y=g(z_1,\dots,z_M)$, where $g(z^*)=\gamma(1,z_1,\dots,z_M)$.\vspace{0.1cm}

Now let us consider the singularities of the hypersurface $V=\{F=0\}$ and its sections. Taking into account the (quasihomogeneous) Euler identity, we get that the closed set $\mathop{\rm Sing}V$ of the singular points of $V$ is given by the system of equations
$$
\frac{\partial F}{\partial x_0}=\dots=\frac{\partial F}{\partial x_M}=\frac{\partial F}{\partial \xi}=0.
$$

{\bf Proposition 1.2.} Let $R\cap\overline{\mathbb P}$ be either the
$\pi_{\mathbb P}$-preimage of a hyperplane in ${\mathbb P}$, or a hypersurface $R_{\gamma}$, defined above. Then
$$
\mathop{\rm dim}\mathop{\rm Sing}(V\cap R)\leqslant\mathop{\rm dim}\mathop{\rm Sing}V+1.
$$.

{\bf Proof.} If $R$ is the $\pi_{\mathbb P}$-preimage of a hyperplane, we repeat word for word the well known argument for the usual projective space. Let $R=R_{\gamma}$. Then the intersection $V\cap R$ is isomorphic to the hypersurface
$$
F(x_0,\dots,x_M,\gamma(x_0,\dots,x_M))=0
$$
in the projective space ${\mathbb P}$, the singularities of which are given by the equations
$$
\frac{\partial F}{\partial x_i}+\frac{\partial F}{\partial \xi}\cdot\frac{\partial\gamma}{\partial x_i}=0,\quad i=0,\dots,M.
$$
Therefore the intersection of $\mathop{\rm Sing}(V\cap R)\subset\overline{\mathbb P}$ with the hypersurface $\left\{\frac{\partial F}{\partial \xi}=0\right\}$ is contained in $\mathop{\rm Sing}V$. Q.E.D. for the proposition.\vspace{0.1cm}

{\bf Definition 1.2.} A non-singular point $o\in V$ is a point {\it of the first type}, if $\frac{\partial F}{\partial \xi}(o)\neq 0$, and {\it of the second type}, if $\frac{\partial F}{\partial \xi}(o)=0$.\vspace{0.1cm}

At a point of the second type the hypersurface, given by the equation
$$
\sum^M_{i=0}x_i\frac{\partial F}{\partial x_i}(o)=0,
$$
is the natural ``tangent hyperplane'' $T_oV\subset\overline{\mathbb P}$. At a point of the first (main) type there is a whole family of candidates for the role of the tangent hyperplane, and they are given by sections of the sheaf ${\cal O}_{\overline{\mathbb P}}(l)$, because they include the coordinate $\xi$: they are of the form $R_{\gamma}$ for suitable polynomials $\gamma$. In order to present the local equation of the hypersurface $V$ more precisely, one needs to consider affine coordinates.\vspace{0.3cm}


{\bf 1.2. Affine coordinates.} Let $o\in V$ be a point. Let us choose projective coordinates on ${\mathbb P}$ in such a way that $o=(1:0:\dots:0:\beta)=(1:0^M:\beta)$ for some $\beta\in{\mathbb C}$. Replacing the coordinate $\xi$ by $\xi'=\xi-\beta x_0^l$, we may assume that $\beta=0$. Now in the affine coordinates $z_1,\dots,z_M,y$ we have $o=(0,\dots,0)$ and $V\cap\{x_0\neq 0\}$ is given by the equation $f=0$, where
$$
f=y^d+a_1(z_*)y^{d-1}+\dots+a_{d-1}(z_*)y+a_d(z_*),
$$
the polynomial $a_i(z_*)$ has degree $il$. Write down
$$
a_i(z_*)=a_{i,0}+a_{i,1}(z_*)+\dots+a_{i,il}(z_*),
$$
where $a_{i,j}$ is a homogeneous polynomial of degree $j$. In particular, $a_{d,0}=0$. The point $o\in V$ is non-singular if and only if the linear form
$$
a_{d-1,0}y+a_{d,1}(z_*)
$$
is not identically zero, and in that case the point $o$ is a point of the first type, if $a_{d-1,0}\neq 0$, and of the second type, if $a_{d-1,0}=0$. If $o\in V$ is a non-singular point of the first type, then $z_1,\dots,z_M$ is a coordinate system on the (affine) tangent space $T_oV$. In particular, every linear subspace $\Lambda\subset T_oV$ is given by a system of linear equations that {\it depend only on} $z_*$, and for any non-zero linear form $h(z_*)$ the intersection
$$
V\cap\{h=0\}
$$
is non-singular at the point $o$. (The last intersection can be understood also as the intersection with the hyperplane $\{h(x_1,\dots,x_M)=0\}$ in $\overline{\mathbb P}$.)\vspace{0.1cm}

Now if $o\in V$ is a non-singular point of the second type, then $z_1,\dots,z_M$, restricted to $T_oV$, are linearly dependent. A typical linear subspace $\Lambda\subset T_oV$ is given by a system of linear equations, one of which is of the form $y-h(z_*)=0$ (where the form $h$ can be identically zero), and the rest depend only on $z_*$.\vspace{0.1cm}

Somewhat abusing the notations, we will use the same symbol $V$ both for the original hypersurface and for its affine part $V\cap\{x_0\neq 0\}$. For a linear form $h(z_*)$ (possibly, identically zero) the symbol $V_h$ means the intersection
$$
V\cap R_{\gamma},
$$
where $\gamma=h(x_1,\dots,x_M)x_0^{l-1}$ in the system of coordinates $(x_0:\dots:x_M:\xi')$, described above. For the affine part of this variety we will use the same symbol $V_h$. Without special comments we consider $V_h$ to be embedded in ${\mathbb P}$ or ${\mathbb A}^M$, depending on the situation.\vspace{0.1cm}

The intersection with a hyperplane $\{h=0\}$ (for a non-zero form $h$) we denote by the symbol of restriction $|_{\{h=0\}}$; again, we use this notation both in the affine and projective context.\vspace{0.3cm}


{\bf 1.3. The regularity conditions.} Let us formulate the regularity conditions for a point $o\in V$. These conditions are slightly different, depending on whether the point $o$ is a non-singular point of the first or second type or a singularity.\vspace{0.1cm}

Assume that the point $o\in V$ is non-singular.\vspace{0.1cm}

(R1) The hypersurface $V_h$ is $\ulcorner(3dl)\slash 8\urcorner$-regular (at the point $o$) for every linear point $h(z_1,\dots,z_M)$.\vspace{0.1cm}

(We use the notations and conventions of Subsection 1.2.)\vspace{0.1cm}

If the point $o\in V$ is non-singular of the second type, that is, $a_{d-1,0}=0$, then, apart from the condition (R1), one more condition is needed. Recall that in that case the tangent hyperplane $T_oV$ is given in the affine coordinates by the equation $a_{d,1}(z_*)=0$.\vspace{0.1cm}

(R1.2) For every linear form $h(z_1,\dots,z_M)$, which is linearly independent with $a_{d,1}(z_*)$, the hypersurface $V_0|_{\{h=0\}}$ is $\ulcorner(3dl)\slash 8\urcorner$-regular.\vspace{0.1cm}

(The last hypersurface is contained in $\{h=0\}={\mathbb A}^{M-1}$.)\vspace{0.1cm}

Now assume that the point $o\in V$ is singular. The affine equation of $V$ at the point $o$ starts with the quadratic form
\begin{equation}\label{02.11.2018.1}
a_{d-2,0}y^2+a_{d-1,1}(z_*)y+a_{d,2}(z_*).
\end{equation}

(R2) The rank of the last form $\geqslant 7$, and the variety $V_0$ is
$\ulcorner dl\slash 2\urcorner$-regular at the point $o$.\vspace{0.1cm}

{\bf Definition 1.1.} The hypersurface $V\in{\cal F}$ is {\it regular}, if at every non-singular point the condition (R1) holds, and at every non-singular point of the second type also the condition (R1.2), and àt every singular point the condition (R2) holds.\vspace{0.1cm}

The set of regular hypersurfaces is denoted by the symbol ${\cal F}_{\rm reg}$. Since every hypersurface $V\in{\cal F}_{\rm reg}$ is either non-singular, or has at most quadratic singularities of rank
$\geqslant 7$, the claim (i) of Theorem 0.1 holds. The claim (ii) of Theorem 0.1 is shown in \S 2.


\section{Codimension of the non-regular set}

The aim of this section is to prove the claim (ii) of Theorem 0.1. In Subsection 2.1 we localize the task: we reduce it to a similar problem for a fixed point $o\in\overline{\mathbb P}$. In Subsection 2.2 we recall the methods of estimating the violations of the regularity condition. In Subsection 2.3 we prove the local estimates, completing the work.\vspace{0.3cm}

{\bf 2.1. The local problem.} Fix a point $o\in\overline{\mathbb P}$,
$o\neq o^*$. Let ${\cal F}(o)\subset{\cal F}$ be the subset (hyperplane) of polynomials that vanish at the point $o$, and
${\cal F}_{\rm reg}(o)\subset{\cal F}(o)$ is the subset of polynomials, satisfying the corresponding regularity condition at that point. Set
$$
{\cal F}_{\rm non-reg}(o)=\overline{{\cal F}(o)\backslash{\cal F}_{\rm reg}(o)}.
$$
Obviously,
$$
{\cal F}\backslash{\cal F}_{\rm reg}\subset\bigcup_{o\in\overline{\mathbb P}\backslash\{o^*\}}{\cal F}_{\rm non-reg}(o),
$$
so that, taking into account the equality $\mathop{\rm codim}({\cal F}(o)\subset{\cal F})=1$, the claim (ii) of Theorem 0.1is implied by the following local fact.\vspace{0.1cm}

{\bf Proposition 2.1.} {\it The following inequality holds:}
$$
\mathop{\rm codim}({\cal F}_{\rm non-reg}(o)\subset{\cal F}(o))\geqslant\frac{(M-4)(M-5)}{2}+M+1.
$$

It is the last inequality that we will show. For each of the reguarity conditions, stated in Subsection 1.3, we have to check that a violation of that condition imposes on the coefficients of the polynomial ${\cal F}\in{\cal F}(o)$ at least $\frac12(M-4)(M-5)+M+1$ independent conditions. As we will see below, the minimal number of independent conditions correspond to the violation of the condition (R2).\vspace{0.3cm}


{\bf 2.2. The methods of estimating the codimension.} We will use two well known methods. The first method was used many times, see \cite[Chapter 3, Section 1]{Pukh13a}. Let, as in Subsection 1.1,
$$
f(z_*)=q_{\mu}(z_*)+\dots+q_N(z_*)
$$
be a polynomial in the affine coordinates $z_1,\dots,z_M$ and $\mu\leqslant k\leqslant M+1-\mu$. By the symbol ${\cal P}_{i,M}$ we denote the space of homogeneous polynomials of degree $i\in{\mathbb Z}_+$ in $z_1,\dots,z_M$, and for $i< j$ by the symbol ${\cal P}_{[i,j],M}$ we denote the product
$$
\prod^j_{a=i}{\cal P}_{a,M},
$$
so that $f\in{\cal P}_{[\mu,N],M}$. Repeating the arguments of Subsection 1.3 in \cite[Chapter 3]{Pukh13a} word for word, we get that the codimension of the set of polynomials $f$ that do not satisfy the condition of  $k$-regularity, with respect to the space ${\cal P}_{[\mu,N],M}$ is not smaller than the number
$$
\mathop{\rm min}\limits_{\mu\leqslant i\leqslant k}
{M-1+\mu\choose i}.
$$

Taking into account the well known behaviour of the binomial coefficients, we conclude that this minimum is realized either for $i=\mu$, or for $i=k$. It is easy to choose, at precisely which value of these two ones it is realized. Note that if it is known that the point $(0,\dots,0)$ is non-singular on the hypersurface $\{f=0\}$, that is, $\mu=1$ and $q_1\not\equiv 0$, then we may fix $q_1$ and restrict $q_i$, $i\geqslant 2$, onto the hyperplane
$\{q_1=0\}$. In that case the codimension is not lower than the number $\left(\begin{array}{c}M\\2\end{array}\right)$ for $k\leqslant M-2$. We will need only this estimate.

The second method is the well known fact that the codimension of the set of quadratic forms of rank $\leqslant r$ in the space of all quadratic forms in $N$ variables is
$$
\frac{(N-r)(N-r+1)}{2}.
$$


{\bf 2.3. Proof of Proposition 2.1.} In the notations of Subsections 1.2, 1.3 let us consider one by one the violations of each regularity conditions at the point $o$. Consider first the hypersurfaces, non-singular at that point.\vspace{0.1cm}

The hypersurface $V=\{f=0\}\subset{\mathbb A}^{M+1}$ is uniquely determined by the set of polynomials
$$
a_1(z_*),\dots,a_d(z_*),
$$
where $\mathop{\rm deg}a_i\leqslant il$. The hypersurface $V_h\subset{\mathbb A}^M$ for some linear form $h(z_*)$ is given by the equation $f_h=0$, where
$$
f_h=h(z_*)^d+a_1(z_*)h(z_*)^{d-1}+\dots+a_d(z_*).
$$
Fix a form $h$, such that $V_h$ is non-singular at the point $o$. Since the polynomial $a_d(z_1,\dots,z_M)$ of degree $\leqslant dl$ is arbitrary, in the presentation
$$
f_h=q_1(z_*)+q_2(z_*)+\dots+q_{dl}(z_*),
$$
where $\mathop{\rm deg}q_i=i$, the homogeneous components $q_i$ are arbitrary and do not depend on each other. Therefore we may apply the first method of estimating the codimension, described in Subsection 2.2. By assumption, $q_1\not\equiv 0$, so that we fix the hyperplane $\{q_1=0\}$ and restrict the polynomial $f_h$ onto this hyperplane:
$$
f_h|_{\{q_1=0\}}=\overline{q}_2+\dots+\overline{q}_{dl}.
$$
Since $\ulcorner(3dl)/8\urcorner\leqslant M-2$, the codimension of the set of polynomials $f_h|_{\{q_1=0\}}$ that do not satisfy the condition of $\ulcorner(3dl)/ 8\urcorner$-regularity, is ${M\choose 2}$. The same is the codimension of the set of non-regular polynomials $f_h$. Since $h$ varies in an $M$-dimensional family, the codimension of the set of polynomials $f$, violating the condition (R1), is
$$
{M\choose 2}-M.
$$
In a similar way we estimate the codimension of the set of polynomials, violating the condition (R1.2): the non-regularity of the hypersurface $V_0|_{\{h=0\}}$ gives
$$
{M-1\choose 2}-M+1
$$
independent conditions for $f$. (The additional codimension +1 comes from the equality $a_{d-1,0}=0$ for a non-singular point of the second type.)\vspace{0.1cm}

Now let us consider the hypersurfaces $V$ which are singular at the point $o$. Since in that case $a_{d-1,0}=0$ and $a_{d,1}(z_*)\equiv 0$, we have $M+1$ additional independent conditions for $f$. If the rank of the quadratic form (\ref{02.11.2018.1}) does not exceed 6, we obtain, therefore,
$$
{M-4\choose 2}+M+1
$$
independent conditions for $f$.\vspace{0.1cm}

Obviously, $\ulcorner dl\slash 2\urcorner\leqslant M-2$. Thus if $V_0$ is a non-regular hypersurface, we obtain
$$
{M+1 \choose 2} + M + 1
$$
independent conditions for $f$.\vspace{0.1cm}

Comparing the results obtained above and choosing the smallest one for $M\geqslant 10$ (which corresponds to the violation of the condition on the rank of the quadratic singularity), we complete the proof of Proposition 2.1 (and the claim (ii) of Theorem 0.1).\vspace{0.1cm}

{\bf Remark 2.1.} We excluded the option $(d,l)=(5,2)$ from consideration for the only reason: it violates the uniformity of the statement of the claim (ii) of Theorem 0.1. For $M=8$ the minimum of the codimension corresponds to the violation of the condition (R1.2) and is equal to 14. With this modification the claims of Theorem 0.1 are true for the values $(d,l)=(5,2)$, too.


\section{Birational superrigidity}

In this section we prove the birational superrigidity of a regular hypersurface $V$. In Subsection 3.1 we recall the key concept of a maximal singularity and exclude certain types of maximal singularities. The remaining types of singularities are classified, after that we start to exclude maximal singularities of ``general position'' (Proposition 3.2). In order to complete this work, we need the technique of hypertangent divisors, which is recalled in Subsection 3.2. Finally, in Subsection 3.3 we exclude maximal singularities of all remaining types, which completes the proof of Theorem 0.1.\vspace{0.3cm}

{\bf 3.1. Maximal singularities.} Let $V\in{\cal F}_{\rm reg}$ be a fixed regular variety. Assume that $V$ is not birationally superrigid. It is well known that in this case on $V$ there is a mobile linear system $\Sigma\subset|nH|$ with a maximal singularity: for some birational morphism $\varphi\colon\widetilde{V}\to V$, where $\widetilde{V}$ is a non-singular projective variety, and some $\varphi$-exceptional prime divisor $Q\subset\widetilde{V}$ the Noether-Fano inequality holds:
$$
\mathop{\rm ord}\nolimits_Q\varphi^*\Sigma>n\cdot a(Q,V),
$$
where $a(Q,V)$ is the discrepancy of $Q$ with respect to $V$ (see, for instance, \cite[Chapter 2]{Pukh13a}). Let $B=\varphi(Q)\subset V$ be the centre of the maximal singularity $Q$. If $B\not\subset\mathop{\rm Sing}V$, then the inequality
\begin{equation}\label{03.11.2018.1}
\mathop{\rm mult}\nolimits_B\Sigma>n
\end{equation}
holds.\vspace{0.1cm}

{\bf Proposition 3.1.} {\it The codimension of the subvariety $B\subset V$ it at least 5.}\vspace{0.1cm}

{\bf Proof.} Assume the converse:
$$
\mathop{\rm codim} (B\subset V)\leqslant 4.
$$
Then $B\not\subset\mathop{\rm Sing}V$, so that the inequality (\ref{03.11.2018.1}) holds. For a general polynomial $g(x_0,\dots,x_M)$ of degree $l$ by Bertini's theorem
$$
\mathop{\rm codim} (\mathop{\rm Sing}\nolimits_g\subset V_g)=
\mathop{\rm codim} (\mathop{\rm Sing}V\subset V)\geqslant 6.
$$
Set $B_g=B\cap V_g$ and $\Sigma_g=\Sigma|_{V_g}$, where $\Sigma_g\subset|nH_g|$ is a mobile linear system on $V_g$, and $H_g$ is the class of a hyperplane section of the hypersurface $V_g\subset{\mathbb P}$. Let $P\subset{\mathbb P}$ be a general 6-plane. The variety
$$
V_P=V_g\cap P
$$
is a non-singular hypersurface in $P\cong{\mathbb P}^6$. On $V_P$ there is a mobile linear system $\Sigma_P=\Sigma_g|_{V_P}\subset |nH_P|$, where $\mathop{\rm Pic}V_P={\mathbb Z}H_P$, and moreover, $\mathop{\rm mult}_{B\cap P}\Sigma_P>n$. However, $B\cap P$ is positive-dimensional, so that the last inequality can not be true (this is a well known fact for any non-singular hypersurface in the projective space, and in fact it is sufficient for the linear system $\Sigma_P$ to be non-empty, see, for instance, \cite[Chapter 2, Section 2]{Pukh13a}). We obtained a contradiction which completes the proof of the proposition.\vspace{0.1cm}

Starting from theis moment, we assume that $\mathop{\rm codim} (B\subset V)\geqslant 5$. Fix a point of general position $o\in B$. There are three options:\vspace{0.1cm}

(1.1) the point $o\not\in\mathop{\rm Sing} V$ is non-singular of the first type,\vspace{0.1cm}

(1.2) the point $o\not\in\mathop{\rm Sing}V$ is non-singular of the second type,\vspace{0.1cm}

(2) the point $o\in\mathop{\rm Sing}V$ is a quadratic singularity.\vspace{0.1cm}

We must exclude each of them.\vspace{0.1cm}

Let us consider first the nonsingular cases. Set $Z=(D_1\circ D_2)$ to be the self-intersection of the linear system $\Sigma$, where $D_1,D_2\in\Sigma$ are general divisors. The effective cycle $Z\sim n^2H^2$ of codimension 2 satisfies the (classical) $4n^2$-inequality
$$
\mathop{\rm mult}\nolimits_oZ>4n^2
$$
and the $8n^2$-inequality
\begin{equation}\label{03.11.2018.2}
\mathop{\rm mult}\nolimits_oZ+\mathop{\rm mult}\nolimits_{\Lambda}Z^+>8n^2,
\end{equation}
here $\Lambda\subset E$ is some linear subspace of codimension 2, $E=\varepsilon^{-1}(o)\subset V^+$ is the exceptional divisor of the blow up
$\varepsilon\colon V^+\to V$ of the point $o$, $E\cong{\mathbb P}^{M-1}$, see, for instance, \cite[Chapter 2, Sections 2,4]{Pukh13a}.\vspace{0.1cm}

{\bf Proposition 3.2.} {\it The case (1.1) does not realize.}\vspace{0.1cm}

{\bf Proof.} Assume the converse: the case 1.1 takes place. In the affine coordinates (see Subsection 1.2) the tangent hyperplane $T_oV$ is given by the equation $a_{d-1,0}y+a_{d,1}(z_*)=0$, where $a_{d-1,0}\neq 0$, so that
$z_1,\dots,z_M$ is a system of coordinates on $T_oV$ and $(z_1:\dots:z_M)$ is a system of homogeneous coordinates on $E={\mathbb P}(T_oV)$. Let
$$
\Lambda=\{h_1(z_*)=h_2(z_*)=0\},
$$
where $h_1,h_2$ are linearly independent forms. Let $h=\lambda_1h_1+\lambda_2h_2$ be a general form in the pencil.\vspace{0.1cm}

Since $\mathop{\rm deg}Z=\mathop{\rm deg}_HZ=dn^2$, the inequality (\ref{03.11.2018.2}) can be re-written in the form
\begin{equation}\label{03.11.2018.3}
\mathop{\rm mult}\nolimits_oZ+\mathop{\rm mult}\nolimits_{\Lambda}Z^+>\frac{8}{d}\mathop{\rm deg}Z.
\end{equation}
This inequality is linear in $Z$. Therefore, there is an irreducible component of the cycle $Z$, satisfying this inequality. In order not to make the notations too complicated, let us simply assume that the cycle $Z$ itself is an irreducible subvariety.\vspace{0.1cm}

{\bf Lemma 3.1.} {\it The subvariety
$$
V\cap\{h_1(z_*)=h_2(z_*)=0\}
$$
is irreducible, non-singular at the point $o$ and not equal to $Z$.}\vspace{0.1cm}

{\bf Remark 3.1.} In the statement of the lemma, the coordinates $z_*$ are considered as {\it affine} coordinates on ${\mathbb A}^{M+1}_{z_*,y}$. We also used $z_*$ above as {\it homogeneous} coordinates on $E$. The subvariety $V\cap\{h_1=h_2=0\}$ is understood as a {\it projective} subvariety in $\overline{\mathbb P}$, that is, the closure of the corresponding affine set. These changes from affine notations to projective ones are obvious and do not require special explanations.\vspace{0.1cm}

{\bf Proof of the lemma.} Non-singularity at the point $o$ is obvious, irreducibility follows from Proposition 1.2 and the assumption on the rank of the quadratic points (the codimension of the singular set $\mathop{\rm Sing}V$ is at least 6). Finally, the subvariety $V\cap\{h_1=h_2=0\}$ has degree $d$ and multiplicity 1 at the point $o$, and its strict transform on $V^+$ has multiplicity precisely 1 along $\Lambda$, so that this subvariety does not satisfy the inequality (\ref{03.11.2018.3}). Therefore, it is not equal to $Z$. Q.E.D. for the lemma.\vspace{0.1cm}

Set
$$
g(z_*)=-\frac{1}{a_{d-1,0}}a_{d,1}(z_*)+h(z_*).
$$
By the lemma, $Z\not\subset V_g$, so that the effective cycle of scheme-theoretic intersection $(Z\circ V_g)$ is well defined. It satisfies the inequality
$$
\mathop{\rm mult}\nolimits_o(Z\circ V_g)>\frac{8}{dl}\mathop{\rm deg}(Z\circ V_g)
$$
(since $V_g\sim lH$ and $V^+_g$ contains $\Lambda$ by the choice of the form
$h$). By the linearity of the last inequality there is an irreducible subvariety $Y\subset V_g$ of codimension 2 (a component of the effective cycle $(Z\circ V_g)$), satisfying the inequality
\begin{equation}\label{05.11.2018.1}
\frac{\mathop{\rm mult}\nolimits_o}{\mathop{\rm deg}}Y>\frac{8}{dl}
\end{equation}
(the symbol $\mathop{\rm mult}\nolimits_o/\mathop{\rm deg}$ means, as usual, the ratio of the multiplicity at the point $o$ to the $H$-degree). The subvariety $Y$ is contained in $V_g$, which is an irreducible hypersurface of degree $dl$ in the projective space ${\mathbb P}$. This hypersurface by the condition (R.1) is $k$-regular at the point $o$, where $k=\ulcorner(3dl)\slash 8\urcorner$.\vspace{0.3cm}


{\bf 3.2. The technique of hypertangent divisors.} We continue our proof of Proposition 3.2. Considering $V_g$ as a hypersurface in the projective space ${\mathbb P}$, let us decompose its equation into components, homogeneous in  $z_*$:
$$
f_g=q_1+q_2+\dots+q_{dl}.
$$
Let us construct the {\it hypertangent linear systems} on $V_g$ at the point $o$:
$$
\Lambda_i=\left\{\sum^i_{j=1}s_{i-j}(q_1+\dots+q_j)|_{V_g}=0\right\},
$$
where $s_a$ independently from each other run through the space ${\cal P}_{a,M}$, see the details and examples in \cite[Chapter 3]{Pukh13a}. By the condition (R1),
$$
\mathop{\rm codim}\nolimits_o(\mathop{\rm Bs}\Lambda_i\subset V_g)=i
$$
for $i=1,\dots,k-1$. Now, applying the technique of hypertangent divisors in the usual way, let us construct a sequence of irreducible subvarieties $Y_2,Y_3,\dots, Y_{k-1}$ of codimension $\mathop{\rm codim}(Y_i\subset V_g)=i$, where $Y_2=Y$ and the last variety in this sequence satisfies the inequality
$$
\frac{\mathop{\rm mult}_o}{\mathop{\rm deg}}Y_{k-1}>\frac{8}{dl}\cdot\frac43\cdot\frac54\cdot
\dots\cdot\frac{k}{k-1}\geqslant 1,
$$
since $8k\geqslant 3dl$ by assumption. This gives the required contradiction and completes the proof of Proposition 3.2.\vspace{0.3cm}


{\bf 3.3. Exclusion of the remaining options.} In order to exclude the cases (1.2) and (2), we repeat the arguments in the proof of Proposition 3.2 with some modifications. We only consider in detail those modifications.\vspace{0.1cm}

{\bf Proposition 3.3.} {\it The case (1.2) does not realize.}\vspace{0.1cm}

{\bf Proof.} Assume the converse. The tangent hyperplane $T_oV$ is given by the equation $a_{d,1}(z_*)=0$. For the subspace $\Lambda\subset E$ there are two options:\vspace{0.1cm}

(1.2.1) $\Lambda$ is given by the equations
$$
y-h_1(z_*)=0,\quad h_2(z_*)=0,
$$
where the forms $h_2$ and $a_{d,1}$ are linearly independent. If $h_1,h_2$ are linearly dependent, then we may assume that $h_1\equiv 0$.\vspace{0.1cm}

(1.2.2) $\Lambda$ is given by the equations
$$
h_1(z_*)=0,\quad h_2(z_*)=0,
$$
where the forms $h_1,h_2$ and $a_{d,1}$ are linearly independent.\vspace{0.1cm}

Assume first that the case (1.2.1) takes place.\vspace{0.1cm}

{\bf Lemma 3.2.}  {\it The subvariety
$$
V\cap\{y-h_1(z_*)=h_2(z_*)=0\}
$$
is irreducible, non-singular at the point $o$ and not equal to $Z$.}\vspace{0.1cm}

{\bf Proof} is completely similar to the proof of Lemma 3.1 and we do not give it.\vspace{0.1cm}

Now set $g(z_*)=h_1(z_*)+\lambda h_2(z_*)$ for a sufficiently general value $\lambda\in{\mathbb C}$. Now the contradiction is obtained by word for word the same arguments as in the case (1.1). We have shown that the option (1.2.1) does not realize.\vspace{0.1cm}

Assume now that the case (1.2.2) takes place. Let
$$
h=\lambda_1h_1+\lambda_2h_2\in\langle h_1,h_2\rangle
$$
be a general form. Again it is easy to check that
$$
Z\not\subset V|_{\{h=0\}},
$$
so that the effective cycle $Z|_{\{h=0\}}=(Z\circ V|_{\{h=0\}})$ has codimension 2 on the irreducible hypersurface $V|_{\{h=0\}}\subset V$ and satisfies the inequality
$$
\mathop{\rm mult}\nolimits_oZ|_{\{h=0\}}>\frac{8}{d}\mathop{\rm deg} Z|_{\{h=0\}}.
$$
By the linearity of this inequality we may assume that the cycle $Z|_{\{h=0\}}$ is an irreducible subvariety. Now if
$$
Z|_{\{h=0\}}\not\subset V_0,
$$
then set $Y$ to be the irreducible component of the effective cycle
$(Z|_{\{h=0\}}\circ V_0)$ with the maximal value of the ratio $\mathop{\rm mult}_o\slash\mathop{\rm deg}$. If, on the contrary, $Z|_{\{h=0\}}\subset V_0$, then set $Y=Z|_{\{h=0\}}$. In any case,
$Y\subset V_0|_{\{h=0\}}$ is a subvariety of codimension 1 or 2, satisfying the inequality (\ref{05.11.2018.1}). Now we argue in the word for word the same way as in the proof of Proposition 3.2, using the condition (R1.2). Q.E.D. for Proposition 3.3.\vspace{0.1cm}

{\bf Proposition 3.4.} {\it The case (2) does not realize.}\vspace{0.1cm}

{\bf Proof.} Assume the converse: $B\subset\mathop{\rm Sing} V$, so that the point $o$ is a quadratic singularity of the variety $V$. By the condition (R2) the point $o\in V$ is a quadratic singularity of rank $\geqslant 7$, so that we can apply the generalized $4n^2$-inequality \cite{Pukh2017a} and conclude that
$$
\mathop{\rm mult}\nolimits_o Z>4n^2
\cdot\mathop{\rm mult}\nolimits_oV=8n^2,
$$
so that
$$
\mathop{\rm mult}\nolimits_oZ>\frac{8}{d}\mathop{\rm deg}Z.
$$
Now we argue as in the proof of Proposition 3.3: by the linearity of the last inequality in $Z$ we may assume that $Z$ is an irreducible subvariety of codimension 2. If $Z\not\subset V_0$, then we set $Y$ to be a component of the effective cycle $(Z\circ V_0)$ with the maximal value of the ratio $\mathop{\rm mult}_o\slash\mathop{\rm deg}$. If $Z\subset V_0$, then we set $Y=Z$. Now we complete the proof in the word for word same way as the proof of Proposition 3.3. The only difference is that now the variety $V_0$ is sibgular at the point $o$: it has a quadratic singularity (of rank $\geqslant 5$), so that if the hypersurface $V_0\subset{\mathbb P}$ is given (in the affine coordinates) by the equation
$$
q_2+q_3+\dots+q_{dl}=0,
$$
then the hypertangent systems are of the form
$$
\Lambda_i=\left\{\sum^i_{j=2}s_{i-j}(q_2+\dots+q_j)|_{V_0}=0\right\},
$$
where $i\geqslant 2$ and the condition for the hypersurface $V_0$ to be $k$-regular, where $k=\ulcorner dl\slash 2\urcorner$, leads in the notations of the proof of Proposition 3.2 to the inequality
$$
\frac{\mathop{\rm mult}_o}{\mathop{\rm deg}}Y_{k-1}>\frac{8}{dl}\cdot\frac54\cdot\frac65
\cdot\dots\frac{k}{k-1}\geqslant 1,
$$
since $2k\geqslant dl$. Q.E.D. for Proposition 3.4.\vspace{0.1cm}

Proof of Theorem 0.1 is complete.\vspace{0.1cm}

{\bf Remark 3.2.} From the very beginning we assumed that $d\geqslant 5$. The double covers are cyclic covers and their superrigidity is well known, see Subsection 0.3. If $d\in\{3,4\}$, then the birational superrigidity follows just from the condition, that the variety $V$ has at most quadratic singularities of rank $\geqslant 7$. Indeed, if $B\subset V$ is the centre of an infinitely near maximal singularity and $\mathop{\rm codim}(B\subset V)\geqslant 3$, then either $B\not\subset\mathop{\rm Sing}V$ and the usual $4n^2$-inequality holds:
$$
\mathop{\rm mult}\nolimits_BZ>4n^2,
$$
or $B\subset\mathop{\rm Sing}V$ and the generalized $4n^2$-inequality holds, which in this case takes the form of the estimate
$$
\mathop{\rm mult}\nolimits_BZ>8n^2.
$$
In any case, $\mathop{\rm mult}_BZ>\mathop{\rm deg}_HZ$, which is impossible
(the linear system $|H|$ is free and defines the finite morphism $\pi\colon V\to{\mathbb P}$). Thus for $d=3$ or 4 the superrigidity holds in essentially weaker assumptions for the variety $V$.


\begin{flushleft}
Department of Mathematical Sciences,\\
The University of Liverpool
\end{flushleft}

\noindent{\it pukh@liverpool.ac.uk}

\end{document}